\documentclass[a4paper]{article}
\usepackage{graphics}
\usepackage{amsmath,amssymb,amsthm,amsfonts}
\usepackage{bbm}
\usepackage{epic,eepic}
\usepackage{psfrag}
\title{Essential edges in Poisson random hypergraphs}

\author{Christina Goldschmidt and James Norris \footnote{Statistical
Laboratory, Centre for Mathematical Sciences, Wilberforce Road,
Cambridge CB3 0WB, UK. \newline Current address for the first author:
Laboratoire de Probabilit\'es et Mod\`eles Al\'eatoires, Universit\'e
Pierre et Marie Curie (Paris 6), 4, place Jussieu, Bo\^\i te courrier
188, 75252 Paris Cedex 05, France; E-mail
\texttt{christina@proba.jussieu.fr} \newline This is a preprint of an
article accepted for publication in \emph{Random Structures and
Algorithms} Copyright \copyright\ (2004) Wiley Periodicals, Inc.} \\
Statistical Laboratory, University of Cambridge}

\newcommand{\E}[1]{\ensuremath{\mathbb{E} \left[#1 \right]}}
\newcommand{\Prob}[1]{\ensuremath{\mathbb{P} \left(#1 \right)}}

\newcommand{\I}[1]{\ensuremath{\mathbbm{1}_{ \{ #1 \} }}}

\newcommand{\Z}{\ensuremath{\mathbb{Z}}}

\newcommand{\ch}[2]{\ensuremath{\left( \begin{smallmatrix} #1 \\ #2
\end{smallmatrix} \right)}}

\newcommand{\sumstack}[2]{\ensuremath{\sum_{\substack{#1 \\ #2}}}}

\newcommand{\sm}{\ensuremath{\setminus}}

\renewcommand{\subset}{\subseteq}
\newcommand{\convdistn}{\ensuremath{\stackrel{d}{\rightarrow}}}
\newcommand{\convprob}{\ensuremath{\stackrel{p}{\rightarrow}}}

\newtheorem{thm}{Theorem}[section]

\newtheorem{lemma}[thm]{Lemma}
\newtheorem{prop}[thm]{Proposition}

\numberwithin{equation}{section}

\begin{document}

\maketitle


\begin{abstract}
\noindent 
Consider a random hypergraph on a set of $N$ vertices in which, for 
$1 \leq k \leq N$, a Poisson$(N \beta_k)$ number of hyperedges is scattered
randomly over all subsets of size $k$. We collapse the hypergraph by running 
the following algorithm to exhaustion: pick a vertex having a $1$-edge and 
remove it; collapse the hyperedges over that vertex onto their remaining
vertices; repeat until there are no $1$-edges left.  We call the
vertices removed in this process \emph{identifiable}.  Also any 
hyperedge all of whose vertices are removed is called \emph{identifiable}.
We say that a hyperedge is \emph{essential} if its removal prior to
collapse would have reduced the number of identifiable vertices.  The
limiting proportions, as $N\rightarrow\infty$, of identifiable
vertices and hyperedges were obtained in \cite{Darling/Norris}. In
this paper, we establish the limiting proportion of essential
hyperedges. We also discuss, in the case of a random graph, the
relation of essential edges to the 2-core of the graph, the maximal
sub-graph with minimal vertex degree $2$.
\end{abstract}



\section{Introduction}

The Poisson random hypergraph model (introduced in
\cite{Darling/Norris}) which is the subject of this paper may be
considered as a step towards developing random combinatorial
structures which one can fit to real-world phenomena.  With this aim
in mind, the class of hypergraphs in which every edge has the same
number of vertices, whilst being a clean and elegant object, may be
too narrow to be useful.  We find it remarkable that, despite the
flexibility afforded by its large number of parameters, the class of
Poisson random hypergraphs admits tractable computations for the
asymptotic size of key structures.

It is conventional in random combinatorics to impose a uniform rather
than a Poisson structure. For large $N$, this makes little difference
so long as one is concerned with ``local" random variables, for
example, the number of edges at a given vertex. A global Poisson
structure is natural probabilistically in that it maximizes
independence. In contrast, in the uniform case, by specifying an exact
total number of hyperedges of a given size, one imposes dependencies
at a global level which may be considered unnatural. For many
questions, including those addressed in this paper, we would expect to
find similar asymptotic behaviour for Poisson and uniform models. This
has already been verified by one of us~\cite{Goldschmidt} in respect
of the numbers of identifiable vertices.

Suppose we have a set of vertices $V$ of size $N$.  Then for
$(\beta_k: k \geq 1)$ a sequence of non-negative real numbers, we
define a Poisson random hypergraph with parameters $(\beta_k: k \geq
1)$ to be a random map $\Lambda: \mathcal{P}(V) \rightarrow \Z^{+}$
such that
\[
\Lambda(A) \sim \mathrm{Poisson} \left(N \beta_k / \ch{N}{k} \right)
\]
whenever $|A| = k$, with $(\Lambda(A),A \subset V)$ independent.  Then
$\Lambda(A)$ is the number of hyperedges over $A$ (so that we are
allowing multiple edges over a set).  We call 1-edges \emph{patches}.
We refer to the case when $\beta_k=0$ for all $k\ge3$ as the graph
case.  Define the generating function $\beta(t) = \sum_{k=1}^{\infty}
\beta_k t^k$.  Throughout this paper, we will assume that
$\sum_{k=1}^{\infty} k\beta_k < \infty$ so that $\beta$ is $C^1$ on
the interval $[0,1]$. Define $t^{*} = \inf \{ t \geq 0: \beta'(t) +
\log(1 - t) < 0 \}$ and note that $t^*\in[0,1)$. We will assume
further that there are no zeros of $\beta'(t) + \log(1 - t)$ in
$[0,t^{*})$.  The case where this last condition fails is explored
further in \cite{Darling/Norris}.

We follow \cite{Darling/Norris}, \cite{Darling/Levin/Norris} in
considering the notion of identifiability.  Any vertex with a patch on
it is identifiable.  Pick such a vertex and delete both the vertex and
the patch.  Collapse all of the other hyperedges over $v$ down onto
their remaining vertices (so that a 3-edge over $\{u,v,w\}$ becomes a
2-edge over $\{u,w\}$, for example).  Continue until there are no more
patches on the hypergraph.  Then the order of this collapse does not
affect the set of vertices eventually removed (see
\cite{Darling/Norris}), called the set of \emph{identifiable}
vertices, which we denote by $V^*(\Lambda)$.  A hyperedge is said to
be identifiable if all of its vertices are identifiable.

\begin{figure}
\begin{center}
\psfrag{$v$}{$v$}
\psfrag{$w$}{$w$}
\resizebox{!}{5cm}{\includegraphics{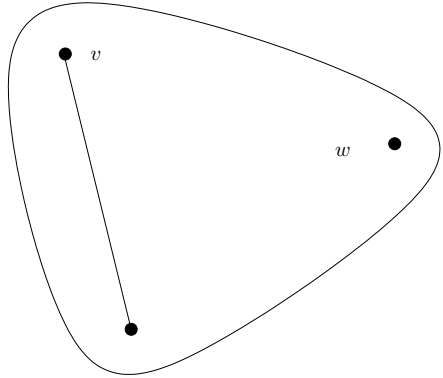}}
\caption{Asymmetry of identifiability}
\label{fig:asymmetry}
\end{center}
\end{figure}

In a hypergraph with no patches there are no identifiable vertices: we
say then that a vertex $w$ is \emph{identifiable from $v$} if it is
identifiable in the hypergraph obtained by adding one patch at $v$.
Note that, except in the graph case, this relation between $w$ and $v$
is not symmetric.  For example, in Figure~\ref{fig:asymmetry}, $w$ is
identifiable from $v$ but $v$ is not identifiable from $w$.  The set
of vertices identifiable from $v$ is called the \emph{domain} of $v$.

We now review some material from \cite{Darling/Norris} and 
\cite{Darling/Levin/Norris} which we will use later.

\begin{thm}\label{thm:darlingnorris}
Let $V_N$ be the number of identifiable vertices and $H_N$ be the
number of identifiable hyperedges in the Poisson random hypergraph on
$N$ vertices.  Then, for all $\epsilon > 0$,
\[
\limsup_{N \rightarrow \infty} N^{-1} \log \Prob{\left| N^{-1} V_N
- t^{*} \right| > \epsilon} < 0
\]
and
\[
\limsup_{N \rightarrow \infty} N^{-1} \log \Prob{\left| N^{-1} H_N
- \beta(t^{*}) + (1 - t^{*}) \log(1 - t^{*}) \right| > \epsilon} < 0.
\]
\end{thm}

Thus, $V_N/N$ and $H_N/N$ have limits in probability which are
attained exponentially fast.

We recall that the Borel($\alpha$) distribution is the distribution of
the total population of a Galton--Watson branching process with
Poisson($\alpha$) offspring distribution.  That is, if $X \sim
\mathrm{Borel}(\alpha)$ then
\begin{align*}
&\Prob{X = n}  = e^{-\alpha n}(\alpha n)^{n-1}/n!, 
                                                \quad n \geq 1 \\
&\Prob{X = \infty}  = p,
\end{align*}
with $p$ the root in $(0,1)$ of $\alpha x + \log(1 - x) = 0$.

\begin{thm}\label{thm:darlinglevinnorris}
Assume that $\beta_1 = 0$.
Let $D_N$ denote the size of the domain of a typical vertex.
Then
\[
D_N \convdistn D
\]
as $N \rightarrow \infty$, where $D$ has the Borel$(2 \beta_2)$
distribution.
\end{thm}

The graph case of this result is well known.  The domain of a vertex
looks like a branching process.  This branching process has almost
surely finite size if $\beta_2 \leq 1/2$ and is infinite with positive
probability if $\beta_2 > 1/2$.


\section{Essential edges}

We say that a hyperedge is \emph{essential} if removing it reduces the
number of identifiable vertices.

Let $\mathcal{E}(N,k)$ be the number of essential $k$-edges in the
Poisson random hypergraph on $N$ vertices.  Let
$\mathcal{E}(N)=\sum_{k=1}^{N} \mathcal{E}(N,k)$ be the total number
of essential edges.  The purpose of this paper is to prove the
following law of large numbers:

\begin{thm} \label{thm:essedges}
As $N\to\infty$, the following limits hold in probability:
\begin{align}
\mathcal{E}(N,k)/N & \to k(1 - t^*)(t^*)^{k-1}\beta_k, \quad k\ge1,
\label{eqn:indconv} \\
\mathcal{E}(N)/N & \rightarrow -(1 - t^{*}) \log(1 -
t^{*}). \label{eqn:fullconv}
\end{align}
\end{thm}

Thus, the limit of $H_N/N$ splits into two parts: $\beta(t^{*})$
corresponding to non-essential identifiable hyperedges and $-(1 -
t^{*}) \log(1 - t^{*}) = (1 - t^{*}) \beta'(t^{*})$ corresponding to
essential hyperedges.


\section{Essential edges in random graphs}

In order to provide intuition about essential edges, we consider first
the case of a random graph with patches. Asymptotically, the
2-edge-structure of the Poisson random graph behaves in the same way as
that of the more-commonly studied binomial model $\mathbb{G}(N, p)$,
with $p = 2\beta_2/N$ (see, for example,
Bollob\'as~\cite{Bollobas} or Janson, \L uczak and
Ruci\'nski~\cite{Janson/Luczak/Rucinski}).
We give a simple calculation for essential patches and then a heuristic
derivation for essential 2-edges based on known results for the 2-core.

In the random graph with patches, a vertex is identifiable if and only
if its component has a patch on some vertex.  So a patch is essential
if and only if it is the only patch on a component.  This enables us
to prove part of our limiting result concerning the expected number of
essential patches in an elementary way. Fix a vertex $v$ and let $D_N$
be the size of the domain of $v$. By
Theorem~\ref{thm:darlinglevinnorris}, $D_N \convdistn D$.  We note
that the limiting Borel$(2 \beta_2)$ distribution for $D$ has
(possibly degenerate) probability generating function $F(s)$ which is
the solution to
\[
F(s) = s \exp(2 \beta_2 (F(s) - 1))
\]
in the range $0 \leq s \leq 1$ (see Harris~\cite{Harris} p.32). 
Let $A$ denote the event that $v$ has a patch and this patch is essential.
Then
\[
\Prob{A} = \E{\Prob{A | D_N}} = \E{ \beta_1 e^{-\beta_1 D_N}} 
\rightarrow \beta_1 F(e^{-\beta_1}).
\]
Now,
\[
F(e^{-\beta_1}) = \exp(-\beta_1 + 2 \beta_2(F(e^{-\beta_1}) - 1)).
\]
But $F(e^{-\beta_1}) = 1 - t^{*}$ is a solution to this equation and
so
\[
\Prob{A} \rightarrow \beta_1 (1 - t^{*}).
\]
Hence, $\E{\mathcal{E}(N,1)/N}$ converges to $\beta_1 (1 - t^{*})$.

A 2-edge in a random graph with patches is essential if and only if
removing it splits its component into two disconnected components,
exactly one of which has a patch on it.  In particular, edges in
cycles are not essential; an edge in a path between two patches is not
essential; edges in components with no patches cannot be essential
(they are not even identifiable).  All edges in a tree-component with
a single patch are essential.

We now describe a connection with the 2-core of a random graph, that is
the maximal subgraph with minimum degree 2.  (The 2-core consists of all
vertices on cycles and on paths between cycles.)  The 2-core of a random
graph is only of size $\Theta(N)$ if the random graph is
super-critical (i.e.~$2 \beta_2 > 1$), so we work with that case.  This
section is intended to provide orientation for our general results in
a context which may be more familiar to some readers.  We do not
attempt to provide a fully rigorous discussion.

It is known (see Pittel~\cite{Pittel}) that the giant component
consists of its own 2-core and a \emph{mantle} of trees, each
sprouting from a different vertex of the 2-core.  The proportion of
vertices in the giant 2-core is a.a.s.~$\theta - 2 \beta_2 \theta (1 -
\theta)$ and the proportion of vertices in the mantle is $2 \beta_2
\theta (1 - \theta)$, where $\theta = \inf\{t \geq 0: 2 \beta_2 t +
\log(1 - t) < 0\}$.  Note that $\theta = t^{*}$ when we have $\beta_1
= 0$ i.e.~no patches.  We can imagine instead that we have an $o(N)$
number of patches, enough to ensure that a.a.s.~one lands on the giant
component.  Then the process of identification basically picks out the
giant component (and a few other smaller tree-like components which we
may neglect because they are of size at most $\mathcal{O}(\log N)$).
Thus, the proportion of identifiable vertices is $t^{*}$ which
corresponds to the size of the giant component.  As edges in cycles
cannot be essential, the essential edges in the giant component must
almost all be found in the mantle.  Because the mantle is like a
forest, the number of edges and vertices in it are approximately equal
and so the number of essential edges scaled by $1/N$ is the same as
the size of the mantle scaled by $1/N$, that is $2 \beta_2 t^*(1 -
t^*)$.


\section{Convergence of expectations}

Let $\Lambda \mathbbm{1}_{\{A\}^C}$ be the Poisson random hypergraph
$\Lambda$ with any hyperedges over the set $A$ removed.
We say that a set $A \subset V$ is \emph{essential for $\Lambda$} if
$\Lambda(A) = 1$ and $V^{*}(\Lambda \mathbbm{1}_{\{A\}^C}) \neq
V^{*}(\Lambda)$.  Thus a set is essential if and only if it has an
essential hyperedge over it.

\begin{prop} \label{prop:equivconds}
A set $A$ is essential if and only if $\Lambda(A) = 1$ and $|A
\sm V^{*}(\Lambda \mathbbm{1}_{\{A\}^C})| = 1$.
\end{prop}

\begin{proof}
It is clear that a set cannot be essential if it has no hyperedges
over it.  Also, it cannot be essential if it has more than
one hyperedge over it.  Recall that the order of deletion does not
affect the set of identifiable edges.  Suppose $|A \sm V^{*}(\Lambda
\mathbbm{1}_{\{A\}^C})| = 0$.  Then everything in $A$ is identifiable
without the hyperedge over $A$ and so if we were to re-introduce the
hyperedge, it would not be essential.  If $|A \sm V^{*}(\Lambda
\mathbbm{1}_{\{A\}^C})| \geq 2$ then replacing the hyperedge over $A$
will not make the two or more vertices identifiable and so $A$ cannot have
been essential.  There remains the case $\Lambda(A) = 1$ with $|A \sm
V^{*}(\Lambda \mathbbm{1}_{\{A\}^C}))| = 1$, when $A$ is obviously essential.
\end{proof}

\begin{lemma} \label{lem:expectation}
For any fixed $k\ge1$, we have $\E{\mathcal{E}(N,k)/N}\to
k(1 - t^*)(t^*)^{k-1}\beta_k$ as $N\to\infty$.
\end{lemma}

\begin{proof}

The numbers of hyperedges on distinct subsets of $V$ are independent
and so $\Lambda \mathbbm{1}_{\{A\}^C}$ has the same law as $\Lambda$
conditioned on $\Lambda(A) = 0$.  Fix $k\ge1$ and choose $A$ with $|A|
= k$.  Let $q^{(1)}(N,k)$ be the probability that $A$ is essential.
Then, by Proposition~\ref{prop:equivconds},
\begin{align*}
q^{(1)}(N,k)
& = \Prob{\Lambda(A) = 1, 
          |A \sm V^{*}(\Lambda \mathbbm{1}_{\{A\}^C})| = 1} \\
& = \Prob{\Lambda(A) = 1} 
    \Prob{|A \sm V^{*}(\Lambda \mathbbm{1}_{\{A\}^C})| = 1} \\
& = \Prob{\Lambda(A) = 1}
    \Prob{|A \sm V^{*}| = 1 | \Lambda(A) = 0} \\
& = \Prob{\Lambda(A) = 1}
    \frac{\Prob{|A \sm V^{*}| = 1}}{\Prob{\Lambda(A) = 0}} \quad
    \text{as $|A \sm V^{*}| = 1$ implies $\Lambda(A) = 0$} \\
& = \frac{N \beta_k}{\ch{N}{k}} \Prob{|A \sm V^{*}| = 1}.
\end{align*}
Hence, 
\begin{equation} \label{eqn:midway}
\E{ \mathcal{E}(N,k)/N}=\ch{N}{k} q^{(1)}(N, k)/N =\beta_k\Prob{|A \sm
V^{*}| = 1}.
\end{equation}

Now, if $|A \sm V^{*}| = 1$, then $A$ must contain $k-1$ identifiable
vertices and one non-identifiable vertex and symmetry implies that,
given $|V^{*}|$, all vertices are equally likely to be identifiable.
So,
\begin{align*}
\Prob{|A \sm V^{*}| = 1 \, \big| \, |V^{*}|} & = \frac{\ch{|V^*|}{k-1}(N -
|V^*|)}{\ch{N}{k}} \\
& = k \left(1 - \frac{|V^{*}|}{N} \right) \frac{|V^*|}{(N - 1)}
\frac{(|V^*| - 1)}{(N - 2)} \cdots
\frac{(|V^*| - k + 2)}{(N - k + 1)} \\
& \convprob k (1 - t^*) (t^*)^{k-1}
\end{align*}
since $|V^*|/N \convprob t^*$.  Because $0 \leq \Prob{|A \sm V^{*}| = 1 \,
\big| \, |V^{*}|} \leq 1$, it follows by bounded convergence that
\begin{equation*} 
\Prob{|A \sm V^*| = 1} = \E{\Prob{|A \sm V^{*}| = 1 \, \big| \, |V^{*}|}}
\rightarrow k(1 - t^*) (t^*)^{k-1}
\end{equation*}
as $N \to \infty$.  Thus, by (\ref{eqn:midway}), $\E{
\mathcal{E}(N,k)/N} \to k(1 - t^*) (t^*)^{k-1} \beta_k$.
\end{proof}


\section{Asymptotic independence}

The key point is to show that the events $\{\text{$A$ is essential}\}$ 
and $\{\text{$B$ is essential}\}$ are asymptotically independent, for distinct
sets $A$ and $B$. Once we have done this, Theorem~\ref{thm:essedges} can be
proved in much the same way as the weak law of large numbers for
independent random variables.

We now proceed as in the proof of Lemma~\ref{lem:expectation}.
Suppose that $|A| = |B| = k$ where $A \neq B$.  Let $q^{(2)}(N,k)$ be
the probability that both $A$ and $B$ are essential.  Then we have
\begin{equation*}
q^{(2)}(N,k) = \Prob{\Lambda(A) = 1, \Lambda(B) = 1, 
                 |A \sm V^{*}(\Lambda\mathbbm{1}_{\{A\}^C})| = 1,
                 |B \sm V^{*}(\Lambda\mathbbm{1}_{\{B\}^C})| = 1}. 
\end{equation*}
Now, if we could deal with $\Lambda\mathbbm{1}_{\{A, B\}^C}$
instead of $\Lambda\mathbbm{1}_{\{A\}^C}$ and
$\Lambda\mathbbm{1}_{\{B\}^C}$ then we would be able to proceed easily
by saying that $\Lambda\mathbbm{1}_{\{A, B\}^C}$ has the same law as
$\Lambda$ conditioned on $\Lambda(A) = 0$ and $\Lambda(B) = 0$.  However,
we must then deal explicitly with the cases where $A$ contributes
towards the identifiability of $B$, or vice versa.  For ease of
notation, write $\tilde{V}^{*}$ for
$V^{*}(\Lambda\mathbbm{1}_{\{A,B\}^C})$.  We commence with some examples.
\begin{figure}
\begin{center}
\psfrag{$A$}{$A$}
\psfrag{$B$}{$B$}
\psfrag{$v$}{$v$}
\psfrag{$w$}{$w$}
\psfrag{H}{$V \sm \tilde{V}^{*}$}
\psfrag{I}{$\tilde{V}^{*}$}
\includegraphics{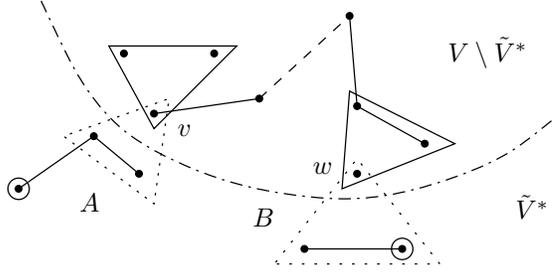}
\caption{Example of two essential edges.  The dot-dashed curve indicates
the boundary between $\tilde{V}^*$ and $V \setminus \tilde{V}^*$.}
\label{fig:simpleess}
\end{center}
\end{figure}
In Figure~\ref{fig:simpleess}, $|A \sm \tilde{V}^{*}| = 1$ and $|B \sm
\tilde{V}^{*}| = 1$.  Both $A$ and $B$ are essential as long as $v$ is
not in the domain of $w$ in $V \sm \tilde{V}^{*}$ and $w$ is not in
the domain of $v$ in $V \sm \tilde{V}^{*}$ (if the dashed edge is
present then $w$ is identifiable from $v$ and so $B$ is not
essential).
\begin{figure}
\begin{center}
\psfrag{$A$}{$A$}
\psfrag{$B$}{$B$}
\psfrag{$v$}{$v$}
\psfrag{$w$}{$w$}
\psfrag{$u$}{$u$}
\psfrag{H}{$V \sm \tilde{V}^{*}$}
\psfrag{I}{$\tilde{V}^{*}$}
\psfrag{J}{$V \setminus V^{*}(\Lambda \mathbbm{1}_{\{B\}^c})$}
\includegraphics{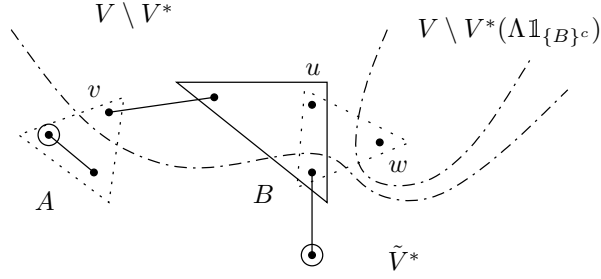}
\caption{Example of two essential edges.  The dot-dashed curves indicate
the boundaries between $\tilde{V}^*$ and $V \setminus \tilde{V}^*$ and
between $V \setminus \tilde{V}^*$ and $V \setminus V^{*}(\Lambda
\mathbbm{1}_{\{B\}^c})$.}
\label{fig:complicatedess}
\end{center} 
\end{figure}
In Figure~\ref{fig:complicatedess}, both $A$ and $B$ are essential but
$|B \sm \tilde{V}^{*}| = 2$.  The important point here is that
precisely one element $u$ of $B \sm \tilde{V}^{*}$ is identifiable
from $v$ and the other is not.

In general, in order to have $|A \sm
V^{*}(\Lambda\mathbbm{1}_{\{A\}^C})| = 1$ and $|B \sm
V^{*}(\Lambda\mathbbm{1}_{\{B\}^C})| = 1$ we must always have either
$|A \sm \tilde{V}^{*}| = 1$ or $|B \sm \tilde{V}^{*}| = 1$ because
otherwise none of $V \sm \tilde{V}^{*}$ would be identifiable when we
reintroduce the edges over $A$ and $B$.  Suppose, without loss of
generality, that $A \sm \tilde{V}^{*} = \{v\}$.  Then we must also
have that all but one of the vertices in $B \sm \tilde{V}^{*}$ are in
the domain of $v$.

\begin{lemma} \label{lem:interaction}
Let $\tilde{\mathcal{D}}(v)$ be the domain of $v$ in the collapsed
hypergraph on $V \sm \tilde{V}^{*}$ and let $\mathcal{D}(v)$ be the domain
of $v$ in the collapsed hypergraph on $V \sm V^{*}$.  Let $\{w\} = B \sm
\tilde{V}^{*}$ when $|B \sm \tilde{V}^{*}| = 1$.  Assume that $v \neq
w$.  Then,
\begin{align}
& \Prob{\Lambda(A) = 1, \Lambda(B) = 1, 
    |A \sm V^{*}(\Lambda\mathbbm{1}_{\{A\}^C})| = 1,
    |B \sm V^{*}(\Lambda\mathbbm{1}_{\{B\}^C})| = 1} \notag \\ 
& = \left(\frac{N \beta_k}{\ch{N}{k}} \right)^2
    \Bigg[ \Prob{|A \sm V^{*}| = 1, |B \sm V^{*}| = 1, 
                  v \not \in \mathcal{D}(w), w \not \in \mathcal{D}(v)} 
						     \notag \\
& \qquad + 2 \sum_{i=2}^{k} 
           \Prob{|A \sm V^{*}| = 1, |B \sm V^{*}| = i, 
                 |\mathcal{D}(v) \cap B \sm V^{*}| = i-1, \Lambda(A) = 0,
                 \Lambda(B) = 0} \Bigg]. \label{eqn:probbothess}
\end{align}
\end{lemma}

\begin{proof}
We have 
\begin{align*}
& \Prob{\Lambda(A) = 1, \Lambda(B) = 1, 
    |A \sm V^{*}(\Lambda\mathbbm{1}_{\{A\}^C})| = 1,
    |B \sm V^{*}(\Lambda\mathbbm{1}_{\{B\}^C})| = 1} \\
& = \Prob{\Lambda(A) = 1, \Lambda(B) = 1, 
    |A \sm \tilde{V}^{*}| = 1, |B \sm \tilde{V}^{*}| = 1, 
    v \not \in \tilde{\mathcal{D}}(w), w \not \in \tilde{\mathcal{D}}(v)} \\
& \qquad + \sum_{i=2}^{k} \Prob{\Lambda(A) = 1, \Lambda(B) = 1, 
    |A \sm \tilde{V}^{*}| = 1, |B \sm \tilde{V}^{*}| = i,
    |\tilde{\mathcal{D}}(v) \cap B \sm \tilde{V}^{*}| = i-1} \\
& \qquad + \sum_{i=2}^{k} \Prob{\Lambda(A) = 1, \Lambda(B) = 1, 
    |B \sm \tilde{V}^{*}| = 1, |A \sm \tilde{V}^{*}| = i,
    |\tilde{\mathcal{D}}(w) \cap A \sm \tilde{V}^{*}| = i-1} \\
& = \Prob{\Lambda(A) = 1} \Prob{\Lambda(B) = 1} \Bigg[
    \Prob{|A \sm \tilde{V}^{*}| = 1, |B \sm \tilde{V}^{*}| = 1, 
    v \not \in \tilde{\mathcal{D}}(w), w \not \in \tilde{\mathcal{D}}(v)}  \\
& \hspace{47mm} 
    + 2 \sum_{i=2}^{k} \Prob{|A \sm \tilde{V}^{*}| = 1,
           |B \sm \tilde{V}^{*}| = i, 
           |\tilde{\mathcal{D}}(v) \cap B \sm \tilde{V}^{*}| = i-1}
    \Bigg] \\
\intertext{}
& = \left( \frac{N \beta_k}{\ch{N}{k}} 
           \exp \left( -N \beta_k/\ch{N}{k} \right) \right)^2 \\
& \qquad \times \Bigg[
         \Prob{|A \sm V^{*}| = 1, |B \sm V^{*}| = 1,
               v \not \in \mathcal{D}(w), 
               w \not \in \mathcal{D}(v) 
               | \Lambda(A) = 0, \Lambda(B) = 0} \\
& \qquad \qquad + 2 \sum_{i=2}^{k} 
         \Prob{|A \sm V^{*}| = 1, |B \sm V^{*}| = i, 
               |\mathcal{D}(v) \cap B \sm V^{*}| = i-1 
               | \Lambda(A) = 0, \Lambda(B) = 0} \Bigg] \\
& = \left( \frac{N \beta_k}{\ch{N}{k}} 
           \exp \left( -N \beta_k/\ch{N}{k} \right) \right)^2 \Bigg[
    \frac{\Prob{|A \sm V^{*}| = 1, |B \sm V^{*}| = 1,
                 v \not \in \mathcal{D}(w), w \not \in \mathcal{D}(v)}}
         {\Prob{\Lambda(A) = 0} \Prob{\Lambda(B) = 0}} \\
& \qquad + 2 \sum_{i=2}^{k} 
    \frac{\Prob{|A \sm V^{*}| = 1, |B \sm V^{*}| = i, 
               |\mathcal{D}(v) \cap B \sm V^{*}| = i-1, \Lambda(A) = 0,
               \Lambda(B) = 0}}
         {\Prob{\Lambda(A) = 0} \Prob{\Lambda(B) = 0}} \Bigg]
\end{align*}
as $|A \sm V^{*}| = 1, |B \sm V^{*}| = 1$ implies $\Lambda(A) = 0,
\Lambda(B) = 0$.  But then cancellation means that the last line is
equal to
\begin{align*}
& \left(\frac{N \beta_k}{\ch{N}{k}} \right)^2
    \Bigg[ \Prob{|A \sm V^{*}| = 1, |B \sm V^{*}| = 1, 
                  v \not \in \mathcal{D}(w), w \not \in \mathcal{D}(v)} \\
& \qquad + 2 \sum_{i=2}^{k} 
           \Prob{|A \sm V^{*}| = 1, |B \sm V^{*}| = i, 
                 |\mathcal{D}(v) \cap B \sm V^{*}| = i-1, \Lambda(A) = 0,
                 \Lambda(B) = 0} \Bigg],
\end{align*}
as required.
\end{proof}
For the moment, we will assume that 
\begin{align}
\Prob{v \not \in \mathcal{D}(w),w \not \in \mathcal{D}(v)}
&\rightarrow 1 \label{eqn:differentdomains} \\
\Prob{|\mathcal{D}(v)\cap B \sm V^{*}| \geq 1} &\rightarrow 0
\label{eqn:intersectdomain} 
\end{align}
as $N \rightarrow \infty$ which, in particular, means that we may
discard the second term in (\ref{eqn:probbothess}).  These results
will be proved later and reflect the fact that the collapse algorithm
will not die out while the vertex domains remain super-critical.
Finally, as we will also show later,
\[
\Prob{|A \sm V^{*}| = 1, |B \sm V^{*}| = 1} \rightarrow (k
(t^{*})^{k-1} (1 - t^{*}))^2
\]
and so
\begin{equation*}
q^{(2)}(N,k) \sim \left( \frac{N \beta_k}{\ch{N}{k}} k (t^{*})^{k-1} 
                                      (1 - t^{*})\right)^2 
\sim q^{(1)}(N,k)^2,
\end{equation*}
where here $\sim$ means that the ratio of the left and right sides
tends to $1$.


\section{Proof of Theorem~\ref{thm:essedges}}

We wish to prove the statements (\ref{eqn:differentdomains}) and
(\ref{eqn:intersectdomain}).  As a first step, we prove

\begin{lemma}
Let $\gamma(t) = (1 - t) \beta''(t)$.  Then
\[
\gamma(t^{*}) \leq 1.
\]
\end{lemma}

\begin{proof}
For $t \in [0,1]$, define $f(t) = 1 - t - \exp(-\beta'(t))$.  Then
\[
f'(t) = \beta''(t) \exp(-\beta'(t)) - 1.
\]
Now $t^{*} = \inf\{t \geq 0: f(t) < 0\}$ and so
$f'(t^{*}) = (1 - t^{*}) \beta''(t^{*}) - 1 \leq 0$.  Hence result.
\end{proof}

As we shall soon see, the 2-edge parameter for the collapsed
hypergraph is approximately $\frac{1}{2} \gamma(t^{*})$ and so this
lemma says that the collapsed hypergraph is nearly sub-critical.

\begin{lemma} \label{lem:probinsamedomain}
Suppose that $v, w \not \in V^{*}$ are chosen uniformly at random.
Recall that $\mathcal{D}(v)$ is the domain of $v$ in the
collapsed hypergraph on $V \sm V^{*}$.  Then
\[
\Prob{w \in \mathcal{D}(v)} \rightarrow 0
\]
as $N \rightarrow \infty$.
\end{lemma}

\begin{proof}
Let $\Lambda^{V^{*}}$ be the hypergraph obtained from $\Lambda$ by
removing all of the vertices in $V^*$, so that $\Lambda^{V^{*}}$ is
the collapsed hypergraph.  Clearly, for $v \in V \sm V^*$,
\begin{equation*}
\Lambda^{V^{*}}(\{v\}) = 0
\end{equation*}
as $v$ is not identifiable.  We need to find the distribution of
$\Lambda^{V^{*}}(A)$ for all $A \in V \sm V^*$.  Suppose that each set
in $\mathcal{P}(V)$ has a corresponding card which gives the number of
hyperedges on it.  Initially, we place the cards face down, so that we
know the sets they represent but not the numbers of hyperedges.
Consider the following slightly different way of looking at the
process of identification.  First turn over the cards corresponding to
all the singleton sets; this tells us which vertices have patches on
them.  Write a list, $\mathcal{L}$ of the vertices with patches.  Now
proceed recursively.  Pick any set with all but one of its vertices in
$\mathcal{L}$ and turn its card over.  Add the last vertex to
$\mathcal{L}$ if there is an edge over the set; if there is no edge,
$\mathcal{L}$ remains unchanged.  Repeat.  The process terminates when
we have run out of sets with all but one of their vertices in
$\mathcal{L}$, so that $\mathcal{L} = V^*$ is the set of identifiable
vertices.  Discard all of the cards which have been turned over, as well as
any corresponding to subsets of $V^*$.  Then there remain only the
cards corresponding to sets with at least two vertices in $V \sm V^*$.
As each carried an independent random variable on its face at the
start and we have not turned any of them over in the process of
identification, the random variables must remain independent of one
another and of what we have seen of the hypergraph on $V^*$.  Thus,
conditional on $V^*$ with $|V^*| = n$, we have that for $A \subset V
\sm V^{*}$ with $|A| = j \geq 2$, 
\[
\Lambda^{V^{*}}(A) = \sumstack{B \supseteq A,}{B \sm V^{*} = A} \Lambda(B)
\sim \mathrm{Poisson} \left( N \sum_{i=0}^{n} \beta_{i+j} \ch{n}{i} /
\ch{N}{i+j} \right)
\]
and these random variables are independent. So $\Lambda^{V^{*}}$ is a
new Poisson random hypergraph on $N - n$ vertices and with parameters
$\beta_k(n,N)$ where
\[
\beta_k(n,N) = \begin{cases}
                 0 & \text{if $k = 1$} \\
                 \frac{N}{N-n} \ch{N-n}{k} \sum_{i=0}^{n} \beta_{i+k}
                 \ch{n}{i} / \ch{N}{i+k} & \text{if $k \geq 2$}.
               \end{cases}
\]
Choose $\rho\in(t^*,1)$.  By Lemma 6.1 of Darling and
Norris~\cite{Darling/Norris}, we obtain
\[
|\beta_2(n,N) - \tfrac{1}{2} \gamma(\tfrac{n}{N})| \leq C(\log N)^2/N
\]
for some constant $C < \infty$ and all $n \leq [\rho N]$.

Furthermore, for each $k$, the number of $k$-edges in the new
hypergraph is certainly bounded by the total number of hyperedges in
the original hypergraph, which had a Poisson($N \beta(1)$)
distribution.  Thus,
\[
\beta_k(n,N) \leq \frac{N \beta(1)}{N - n} \leq \frac{\beta(1)}{1 -
\rho}
\]
for all $k \geq 2$ and all $n \leq [\rho N]$.

Now let $C_t$ have the distribution of the size of the domain
of a vertex in a Poisson random hypergraph with parameters $\beta_1 = 0$,
$\beta_2 = \frac{1}{2} \gamma(t) + C(\log N)^2/N$ and $\beta_k =
\beta(1)/(1 - \rho)$ for $k \geq 3$.  Then by
Theorem~\ref{thm:darlinglevinnorris} and an obvious comparison argument,
\[
C_t \convdistn \mathrm{Borel}(\gamma(t))
\]
as $N \rightarrow \infty$.  Also, if $|V^{*}| = Nt$ for $t < \rho$
then $|\mathcal{D}(v)|$ is stochastically dominated by $C_t$.

For any $\delta > 0$, choose $\epsilon > 0$ small enough and $R
< \infty$ large enough that, firstly, $t^* + \epsilon < \rho$ and,
secondly, that
\[
\Prob{\mathrm{Borel}(\gamma(t^{*} + \epsilon)) \leq R} \geq 1 -
  \delta
\]
(we can do this because $\gamma(t^{*}) \leq 1$ and $\gamma$ is continuous).
Hence,
\begin{align*}
\Prob{w \in \mathcal{D}(v)}
& = \Prob{w \in \mathcal{D}(v), \left| \frac{|V^{*}|}{N} - t^{*} \right| 
          \leq \epsilon}
    + \Prob{w \in \mathcal{D}(v), \left| \frac{|V^{*}|}{N} - t^{*} \right| 
            > \epsilon} \\
& \leq \Prob{w \in \mathcal{D}(v), |\mathcal{D}(v)| \leq R, 
             \left| \frac{|V^{*}|}{N} - t^{*} \right| \leq \epsilon} \\
& \qquad + \Prob{|\mathcal{D}(v)| > R, 
               \left| \frac{|V^{*}|}{N} - t^{*} \right| \leq \epsilon}
       + \Prob{\left| \frac{|V^{*}|}{N} - t^{*} \right| > \epsilon} \\
& \leq \frac{R}{N(1 - t^{*} - \epsilon)} + \Prob{C_{t^{*} + \epsilon}
> R} + \Prob{\left| \frac{|V^{*}|}{N} - t^{*} \right| > \epsilon} \\
& \rightarrow \Prob{\mathrm{Borel}(\gamma(t^{*} + \epsilon)) > R}
\quad \text{as $N \rightarrow \infty$}
\end{align*}
But this last quantity is less than $\delta$ and so we are done.
\end{proof}

Finally, we give a technical lemma.

\begin{lemma} \label{lem:sumconv}
Suppose that for $k \geq 1$ and $N \geq 1$, $X_{k,N}$ are non-negative
random variables satisfying
\begin{equation} \label{eqn:basicconv}
\frac{1}{N} X_{k,N} \convprob x_k
\end{equation}
as $N \to \infty$, for all $k$, where
\begin{equation} \label{eqn:xkbound}
\sum_{k=1}^{\infty} x_k < \infty.
\end{equation}
Suppose in addition that for each $k$ there exists $y_k$ such that
\begin{equation} \label{eqn:expbound}
\E{\frac{1}{N} X_{k,N}} \leq y_k
\end{equation}
for all $N$ with
\begin{equation} \label{eqn:ykbound}
\sum_{k=1}^{\infty} y_k < \infty.
\end{equation}
Then
\begin{equation*}
\frac{1}{N} \sum_{k=1}^{N} X_{k,N} \convprob \sum_{k=1}^{\infty} x_k
\end{equation*}
as $N \to \infty$.
\end{lemma}

\begin{proof}
Let $\epsilon, \delta > 0$.  By (\ref{eqn:xkbound}) and
(\ref{eqn:ykbound}), we can find $k_0$ sufficiently large that
\begin{equation} \label{eqn:Markovpart}
\sum_{k=k_0 + 1}^{\infty} x_k + \sum_{k=k_0 + 1}^{\infty} y_k < \delta
\epsilon / 4.
\end{equation}
Moreover, as $k_0$ is fixed and finite, by (\ref{eqn:basicconv})we
have that $N^{-1} \sum_{k=1}^{k_0} X_{k,N} \convprob \sum_{k=1}^{k_0}
x_k$ and so, for $N$ sufficiently large,
\begin{equation} \label{eqn:finitepart}
\Prob{ \left| \frac{1}{N} \sum_{k=1}^{k_0} X_{k,N} - \sum_{k=1}^{k_0}
x_k \right| > \epsilon/2} < \delta/2.
\end{equation}
Now,
\begin{align*}
& \Prob{\left| \frac{1}{N} \sum_{k=1}^{N} X_{k,N} - \sum_{k=1}^{\infty} x_k
\right| > \epsilon} \\
& \qquad \leq \Prob{ \left| \frac{1}{N} \sum_{k=1}^{k_0} X_{k,N}
- \sum_{k=1}^{k_0} x_k \right| > \epsilon/2} 
+ \Prob{\left| \frac{1}{N} \sum_{k=k_0+1}^{N} X_{k,N} -
\sum_{k=k_0+1}^{\infty} x_k \right| > \epsilon/2} \\
& \qquad \leq \Prob{ \left| \frac{1}{N} \sum_{k=1}^{k_0} X_{k,N} -
\sum_{k=1}^{k_0} x_k \right| > \epsilon/2} 
+ \frac{2}{\epsilon} \left\{ \sum_{k=k_0 + 1}^{\infty} y_k +
\sum_{k=k_0 + 1}^{\infty} x_k \right\} 
\end{align*}
by Markov's inequality and (\ref{eqn:expbound}).  But by
(\ref{eqn:Markovpart}) and (\ref{eqn:finitepart}) this last expression
is less than $\delta$ for all $N$ sufficiently large.
\end{proof}

Now we are ready to prove Theorem~\ref{thm:essedges}.

\emph{Proof of Theorem~\ref{thm:essedges}.}  Firstly observe that
(\ref{eqn:fullconv}) follows from (\ref{eqn:indconv}) by applying
Lemma~\ref{lem:sumconv} with $X_{k,N} = \mathcal{E}(N,k)$, $x_k = k(1
- t^*)(t^*)^{k-1} \beta_k$ and $y_k = \beta_k$ (the number of
essential $k$-edges per vertex is bounded in expectation by the total
number of $k$-edges per vertex), where we have $\sum_{k=1}^{\infty}k(1
- t^*)(t^*)^{k-1} \beta_k < \infty$ and $\sum_{k=1}^{\infty}\beta_k
< \infty$ by assumption.

It remains to prove (\ref{eqn:indconv}).  Take $\epsilon > 0$ and fix
$k \geq 1$.  Note that
\begin{equation*}
\frac{1}{N} \mathcal{E}(N,k) = \frac{1}{N}
\sumstack{A \subset V,}{|A| = k} \I{\text{$A$ is essential}}.
\end{equation*}
By Chebyshev's inequality, for $N$ sufficiently large,
\begin{align*}
& \Prob{\Bigg| \frac{1}{N}  
        \sumstack{A \subset V,}{|A| = k} \I{\text{$A$ is essential}} 
        - (1 - t^{*}) k \beta_k (t^{*})^{k-1} \Bigg| > \epsilon} \\
& \leq \frac{1}{\epsilon^2} \E{\left(\frac{1}{N} 
        \sumstack{A \subset V,}{|A| = k} \I{\text{$A$ is essential}} 
        - (1 - t^{*}) k \beta_k (t^{*})^{k-1} \right)^2} \\
& = \frac{1}{\epsilon^2} 
     \left\{\frac{1}{N^2}  
            \sumstack{A,B \subset V, A \neq B}{|A|=|B|=k} 
            \Prob{\text{$A$ is essential, $B$ is essential}}
            - ((1 - t^{*}) k \beta_k (t^{*})^{k - 1})^2 
    \right\} \\
& \qquad + \frac{1}{\epsilon^2} 
    \left\{ \frac{1}{N^2} \E{\mathcal{E}(N,k)} 
            + 2(1 - t^{*}) k \beta_k (t^{*})^{k-1}
              \left[ (1 - t^{*}) k \beta_k (t^{*})^{k-1} 
                     - \frac{1}{N} \E{\mathcal{E}(N,k)}
              \right] 
    \right\},
\end{align*}
where the second term in braces tends to 0 as $N \rightarrow \infty$,
by Lemma~\ref{lem:expectation}.  By Lemma~\ref{lem:interaction}, 
\begin{align*}
& \frac{1}{N^2} \sumstack{A,B \subset V, A \neq B}{|A|=|B|=k} 
                \Prob{\text{$A$ is essential, $B$ is essential}} \\
& = \frac{1}{N^2} \sumstack{A,B \subset V, A \neq B}{|A|=|B|=k}
            \frac{N^2 \beta_k^2}{\ch{N}{k}^2}
            \Bigg\{
            \Prob{|A \sm V^{*}| = 1, |B \sm V^{*}| = 1, v \not \in
            \mathcal{D}(w), w \not \in \mathcal{D}(v)} \\
& \hspace{2cm} + 2 \sum_{i=2}^{k} 
           \Prob{|A \sm V^{*}| = 1, |B \sm V^{*}| = i, 
                 |\mathcal{D}(v) \cap B \sm V^{*}| = i-1, \Lambda(A) = 0,
                 \Lambda(B) = 0} \Bigg\}.
\end{align*}
Now, by Lemma~\ref{lem:probinsamedomain}, 
$\Prob{v \not \in \mathcal{D}(w), w \not \in \mathcal{D}(v)}
\rightarrow 1$ and $\Prob{|\mathcal{D}(v) \cap B \sm V^{*}| \geq 1}
\rightarrow 0$ as $N \rightarrow \infty$.  Thus, we are really interested in
\begin{align}
& \frac{\beta_k^2}{\ch{N}{k}^2}
    \sumstack{A,B \subset V, A \neq B}{|A|=|B|=k}
    \Prob{|A \sm V^{*}| = 1, |B \sm V^{*}| = 1} \notag \\
& \qquad = \frac{\beta_k^2}{\ch{N}{k}^2} 
    \sumstack{A,B \subset V, A \cap B = \emptyset}{|A|=|B|=k}
    \Prob{|A \sm V^{*}| = 1, |B \sm V^{*}| = 1} \notag \\
& \qquad \qquad + \frac{\beta_k^2}{\ch{N}{k}^2} 
    \sum_{i=1}^{k-1} \sumstack{A,B \subset V, |A \cap B| = i}{|A|=|B|=k}
    \Prob{|A \sm V^{*}| = 1, |B \sm V^{*}| = 1}. \label{eqn:intersections}
\end{align}
In the case where $A$ and $B$ are disjoint, by a similar argument to
that used in the proof of Lemma~\ref{lem:expectation}, we have
\begin{equation*}
\Prob{|A \sm V^*| = 1, |B \sm V^{*}| = 1} = \E{\frac{(N - |V^*|)(N -
|V^*| - 1) \ch{|V^*|}{k-1} \ch{|V^*| - k + 1}{k - 1}}{\ch{N}{k}
 \ch{N-k}{k}}}.
\end{equation*}
If $|A \cap B| = i$ for $1 \leq i \leq k-1$, then the two
non-identifiable vertices must lie outside the intersection (otherwise
they are the same vertex and so $v = w$).  Thus,
\begin{align*}
& \Prob{|A \sm V^*| = 1, |B \sm V^{*}| = 1} \\
& \qquad \qquad = \E{\frac{(N -
|V^*|)\ch{|V^*|}{k - i - 1} \ch{|V^*| - k + i + 1}{i} (N - |V^*| - 1)
 \ch{|V^*| - k + 1}{k - i - 1}}{\ch{N}{k} \ch{k}{i} \ch{N - k}{k - i}}}.
\end{align*}
So (\ref{eqn:intersections}) is equal to
\begin{align*}
& \frac{\beta_k^2}{\ch{N}{k}^2} \E{(N - |V^*|)(N -
|V^*| - 1) \ch{|V^*|}{k-1} \ch{|V^*| - k + 1}{k - 1}} \\
& \qquad + \frac{\beta_k^2}{\ch{N}{k}^2} \sum_{i=1}^{k-1} \E{(N -
|V^*|)\ch{|V^*|}{k - i - 1} \ch{|V^*| - k + i + 1}{i} (N - |V^*| - 1)
 \ch{|V^*| - k + 1}{k - i - 1}}.
\end{align*}
As $\frac{|V^{*}|}{N} \rightarrow t^{*}$ in probability, the first of
these two terms converges to
\[
(k \beta_k (1 - t^{*}) (t^{*})^{k-1})^2,
\]
by bounded convergence.  It remains to show that the second term
converges to 0.  We have
\begin{align*}
& \sum_{i=1}^{k-1} \E{\frac{(N - |V^*|)\ch{|V^*|}{k - i - 1}
                 \ch{|V^*| - k + i + 1}{i} (N - |V^*| - 1)  \ch{|V^*|
                 - k + 1}{k - i - 1}}{\ch{N}{k}^2}} \\
& \quad \leq \sum_{i=1}^{k-1} \frac{(k!)^2}{i! ((k-i-1)!)^2 (N - k + i
                 + 1)!} \\
& \quad \rightarrow 0
\end{align*}
as $N \rightarrow \infty$.  The result follows.
\hfill $\Box$


\end{document}